\documentclass[amsthm]{elsart}

\usepackage{amsmath,url}

\usepackage[OT2,T1]{fontenc}
\DeclareSymbolFont{cyrletters}{OT2}{wncyr}{m}{n}
\DeclareMathSymbol{\shuffle}{\mathbin}{cyrletters}{"78}

\newcommand{\Ainf}{$A_\infty$}
\newcommand{\Cinf}{$C_\infty$}

\newcommand{\Om}{\Omega}
\newcommand{\om}{\omega}
\newcommand{\eps}{\varepsilon}
\newcommand{\bull}{\bullet}
\newcommand{\p}{\partial}
\DeclareMathOperator{\N}{NC}

\renewcommand{\o}{\otimes}

\DeclareMathOperator{\Tot}{Tot}
\DeclareMathOperator{\TotTW}{\textup{Tot}_{\textup{TW}}}

\newcommand{\tint}{{\textstyle\int}}

\renewcommand{\k}{\mathbf{K}}

\newcommand{\T}{\mathcal{T}}

\newcommand{\id}{1}

\begin{document}

\begin{frontmatter}

\title{Transferring homotopy commutative algebraic structures}

\author[USTC]{Xue Zhi Cheng},
\ead{xuezhi.cheng@gmail.com}
\author[NU]{Ezra Getzler\thanksref{nsf}}
\ead{getzler@northwestern.edu}

\address[USTC]{Department of Mathematics, USTC, Hefei 230026, China, \\ 
  and \\ CMTP, USTC-SIAS, Shanghai 201315, China}

\address[NU]{Department of Mathematics, Northwestern University,
  Evanston, IL 60208, USA}

\begin{abstract}
  We show that the sum over planar trees formula of Kontsevich and
  Soibelman transfers \Cinf-structures along a contraction.  Applying
  this result to a cosimplicial commutative algebra $A^\bull$ over a
  field of characteristic zero, we exhibit a canonical \Cinf-structure
  on $\Tot(A^\bull)$, which is unital if $A^\bull$ is; in particular,
  we obtain a canonical \Cinf-structure on the cochain complex of a
  simplicial set.
\end{abstract}

\thanks[nsf]{Partially supported by NSF Grants DMS-0072508 and
  DMS-0505669.}

\end{frontmatter}

\thispagestyle{empty}

\subsection*{Contractions}

Consider cochain complexes $A$ and $B$ over a field $\k$ of
characteristic $0$.
\begin{defn}
  A contraction $(f,g,H)$ from $A$ to $B$ consists of chain maps
  $f:A\to B$ and $g:B\to A$ and a chain homotopy $H:A\to A[1]$, such
  $f\circ g = \id_B$, and
  \begin{equation}
    g\circ f - \id_A = d_AH + Hd_A .
  \end{equation}
\end{defn}

If $A$ is an \Ainf-algebra, there is a natural way to transfer an
\Ainf-structure along a contraction: the induced \Ainf-structure on
$B$ is given by an explicit formula due to Kontsevich and
Soibelman~\cite{KS}. (The case in which $A$ is a differential graded
algebra and the differential on $B$ vanishes is due to Kadeishvili
\cite{Kadeishvili}.)

\subsection*{\Ainf-algebras and \Ainf-morphisms}

The following definition of an \Ainf-structure on a cochain complex
$A$ is not the original one (Stasheff \cite{Stasheff}), which was
formulated on the shifted complex $A[1]$. The grading convention
adopted here has the advantage of supressing many signs. For a
discussion of these signs, see Getzler and Jones~\cite{GJ}.
\begin{defn}
  An \Ainf-algebra is a cochain complex $A$ over a field $\k$,
  equipped with multilinear mappings
  \begin{equation*}
    m_n^A : A^{\o n} \to A , \quad n\ge1 ,
  \end{equation*}
  of degree $1$, such $d=m_1^A$ and for all $n\ge 1$,
  \begin{equation}
   \sum_{k=1}^{n}\sum_{j=0}^{n-k} m_{n-k+1}^{A}\circ(1_A^{\o j}\o
      m_{k}^{A}\o 1_A^{\o n-j-k})=0 .
      \tag{$\text{I}_n$}
  \end{equation}
\end{defn}

The first of the above equations states that the differential $m_1^A$
has square zero,
\begin{equation}
  m_1^A \circ m_1^A=0 ,
  \tag{$\text{I}_1$}
\end{equation}
the second, that the bilinear product $m_2^A$ on $A$ has $m_1^A$ as a
derivation,
\begin{equation}
  m_1^A \circ m_2^A + m_2^A\circ(m_1^A\o 1_A+ 1_A\o m_1^A) = 0 ,
  \tag{$\text{I}_2$}
\end{equation}
and the third, that $m_2^A$ is associative up to the homotopy $m_3^A$,
\begin{multline}
  m_2^A\circ(1_A\o m_2^A+m_2^A\o1_A) + m_1^A\circ m_3^A \\
  + m_3^A\circ(m_1^A\o1_A^{\o 2}+1_A\o m_1^A\o1_A+1_A^{\o2}\o m_1^A) = 0 .
  \tag{$\text{I}_3$}
\end{multline}
If $A$ is an \Ainf-algebra such that $m_n^A=0$ for $n>2$, then $A[1]$
is a differential graded associative algebra, with product
\begin{equation}
  \label{ab}
  ab = (-1)^{|a|+1} m_2^A(a,b) .
\end{equation}
Here, as elsewhere in this paper, $|a|$ denotes the degree $k$ of
the homogeneous element $a\in A^k$.

\begin{defn}
  An \Ainf-morphism $F:A\to B$ of \Ainf-algebras is a series of
  multi-linear maps $F_n:A^{\o n}\to B$, $n\ge1$, of degree $0$ such
  that for all $n\ge 1$
  \begin{multline}
    \sum_{k=1}^n \sum_{\substack{n_1+\dots+n_k=n \\ n_i>0}}
    m_k^B(F_{n_1}\o\dots\o F_{n_{k}}) \\
    = \sum_{k=1}^n \sum_{j=0}^{n-k} F_{n-k+1}(1_A^{\o j}\o
    m_k^A\o1_A^{\o n-j-k}) .  \tag{$\text{II}_n$}
  \end{multline}

  An \Ainf-morphism $F$ is called a quasi-isomorphism if $F_1$ is a
  quasi-isomorphism.
\end{defn}

The first of the above equations states that $F_1$ is a cochain map,
\begin{equation}
  m_1^B \circ F_1 = F_1 \circ m_1^A ,
  \tag{$\text{II}_1$}
\end{equation}
and the second, that $F_1$ is a morphism for the products $m_2^A$ and
$m_2^B$ up to the homotopy $F_2$,
\begin{equation}
  m_2^B\circ(F_1\o F_1) + m_1^B \circ F_2 = F_2\circ(m_1^A\o1_A+1_A\o
  m_1^A) + F_1\circ m_2^A .
  \tag{$\text{II}_2$}
\end{equation}

\subsection*{Rooted planar trees and the transfer formula}

The formula of Kontsevich and Soibelman for the transfer of an
\Ainf-structure is expressed using the language of planar trees.
\begin{defn}
  A rooted tree $T$ is a directed, simply connected graph, with
  vertices $V(T)$ and edges $E(T)$, such that each vertex is the
  source of at most one edge.
\end{defn}

Denote the source of an edge by $\p_s(e)$, and its target by
$\p_t(e)$.

\begin{defn}
  A planar tree is a tree with a total order on the set $\p_t^{-1}(v)$
  for all $v\in V(T)$.
\end{defn}

The valence $|v|=|\p_t^{-1}(v)|$ of a vertex $v$ is the number of
edges having $v$ as target. A tail is a vertex of valence $0$. There
is a unique vertex which is the source of no edge, called the root of
the tree. An edge $e\in E(T)$ is interior if it does not meet a tail
or the root.

Let $\T_n$ be the set of isomorphism classes of planar trees with $n$
tails, each vertex having valence greater than $1$.

Let $A$ be an \Ainf-algebra, and let $(f,g,H)$ be a contraction from
$A$ to a cochain complex $B$. To each planar tree $T$ in $\T_n$,
define an operation $m_T: B^{\o n}\to B[1]$: Assign to each tail of
$T$ the map $g$, to each vertex $v$ of valence $k$ the map $m^A_k$, to
each interior edge the homotopy $H$, and to the root the map $f$.
Moving down the tree from the tails to the root, one reads off the map
$m_T$ as the composition of these different assignments.

\begin{figure}[hbp]
  \centering
  \begin{picture}(175,180)(100,0)
    \put(190,111){\circle*{4}}
    \put(195,105){$m_{4}$}
    \put(145,155){\circle*{4}}
    \put(175,155){\circle*{4}}
    \put(205,155){\circle*{4}}
    \put(235,155){\circle*{4}}
    \put(190,65){\line(0,1){45}}
    \put(190,110){\line(-1,1){45}}
    \put(190,110){\line(1,1){45}}
    \put(190,110){\line(1,3){15}}
    \put(190,110){\line(-1,3){15}}
    \put(141,162){$g$}
    \put(171,162){$g$}
    \put(202,162){$g$}
    \put(232,162){$g$}
    \put(190,65){\circle*{4}}
    \put(195,60){$m_{3}$}
    \put(145,110){\circle*{4}}
    \put(235,110){\circle*{4}}
    \put(190,65){\line(0,-1){45}}
    \put(190,65){\line(-1,1){45}}
    \put(190,65){\line(1,1){45}}
    \put(141,117){$g$}
    \put(232,117){$g$}
    \put(190,20){\circle*{4}}
    \put(186,6){$f$}
  \end{picture}  
  \caption{The tree $T$}
  \label{T}
\end{figure}
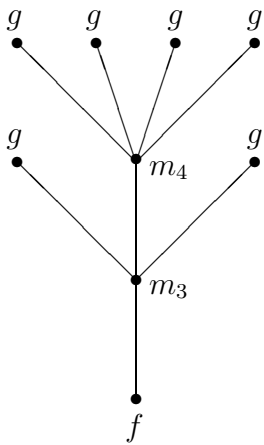

For example, the tree $T$ of Fig.~\ref{T} yields the operation
\begin{equation*}
  m_T = f\circ m_{3}\circ(g\o (H\circ m_{4}\circ (g\o g\o g\o g))\o g)
  .
\end{equation*}

Define operations $G_T: B^{\o n}\to A$ in the same way, except that
now the operator $H$ is assigned to the root instead of $f$.

We may now state the theorem of Kontsevich and Soibelman \cite{KS}.
\begin{thm}
  \label{KS}
  The sequence of operations $m_1^B=d_B$ and
  \begin{equation*}
    m_n^B = \sum_{T\in\T_n} m_T , \quad n>1 ,
  \end{equation*}
  are an \Ainf-structure on $B$.
  
  The sequence of operations $G_1=g$ and
  \begin{equation*}
    G_n = \sum_{T\in\T_n} G_T , \quad n>1 ,
  \end{equation*}
  are an \Ainf-morphism from $B$ to $A$.
\end{thm}

In proving theorems about the transferred \Ainf-structure on $B$, we
will make use of the following inductive formula for $m_n^B$ and
$G_n$, whose proof is clear.
\begin{lem}
  \label{inductive}
  For $n>1$, we have
  \begin{align*}
    m_n^B &= \sum_{k=2}^n \sum_{\substack{n_1+\dots+n_k=n \\ n_i>0}} f
    \circ m_k^A \circ ( G_{n_1} \o \dots \o G_{n_k} ) \intertext{and}
    G_n &= \sum_{k=2}^n \sum_{\substack{n_1+\dots+n_k=n \\ n_i>0}} H
    \circ m_k^A \circ ( G_{n_1} \o \dots \o G_{n_k} ) .
  \end{align*}
\end{lem}
  
\subsection*{Transfer of unital \Ainf-structures}

\begin{defn}
  An identity for an \Ainf-algebra $A$ is an element $e_A\in A$, of
  degree $-1$, such that $m_1^A(e_A)=0$,
  $m_2^A(e_A,a)=(-1)^{|a|+1}m_2^A(a,e_A)=a$, and
  \begin{equation*}
    m_n^A(a_1,\dots,e_A,\dots,a_n) = 0 , \quad n>2 .
  \end{equation*}
  A unital \Ainf-algebra is an \Ainf-algebra $A$ with identity $e_A$.
\end{defn}

\begin{defn}
  An \Ainf-morphism $G$ of unital \Ainf-algebras from $A$ to $B$ is
  unital if $G_1(e_A)=e_B$ and
  \begin{equation*}
    G_n(a_1,\dots,e_A,\dots,a_n) = 0 , \quad n>1 .
  \end{equation*}
\end{defn}

The following is a consequence of the explicit formulas of
Theorem~\ref{KS}.
\begin{thm}
  \label{unital}
  Let $A$ be a unital \Ainf-algebra, and let $(f,g,H)$ be a
  contraction from $A$ to the cochain complex $B$. If the contraction
  satisfies the side conditions
  \begin{align*}
    & f\circ H = 0 , & & H\circ H = 0 , & & H(e_A) = 0 ,
  \end{align*}
  then the induced \Ainf-structure on $B$ is unital, with
  identity $e_B=f(e_A)$, and $G$ is a unital \Ainf-morphism.
\end{thm}
\begin{proof}
  First, observe that these side conditions imply that $H\circ g=0$:
  indeed, we have
  \begin{equation*}
    [g\circ f,H] = [\id_A-Hd_A-d_AH,H] = 0 ,
  \end{equation*}
  and hence
  \begin{equation*}
    H\circ g = H\circ g\circ f\circ g = g\circ f\circ H\circ g = 0 .
  \end{equation*}

  We prove the theorem by induction on $n$. For $n=1$, we have
  \begin{equation*}
    G_1(e_B) = g(e_B) = g\circ f(e_A) = ( \id_A - d_A H - H d_A ) e_A
    = e_A ,
  \end{equation*}
  and
  \begin{equation*}
    m_1^B(e_B) = f\circ m_1^A\circ g(e_B) = f\circ m_1^A(e_A) = 0 .
  \end{equation*}
  
  For $n>1$, we have by Lemma \ref{inductive}
  \begin{multline*}
    G_n^B(b_1,\dots,b_n) \\
    = \sum_{k=2}^n \sum_{\substack{n_1+\dots+n_k=n \\ n_i>0}} f \circ
    m_k^A( G_{n_1}(b_1,\dots,b_{n_1}),\dots,G_{n_k}(b_{n-n_k+1},\dots,b_n)) .
  \end{multline*}
  If $b_j=e_B$, a given term on the right-hand side vanishes, by the
  induction hypothesis, unless $j=n_1+\dots+n_i$ and $n_i=1$. In this
  case, the $i$th argument of $m_k^A$ is $G_1(e_B)=e_A$, and hence the
  term vanishes unless $k=2$, and either $j=1$ or $j=n$. In the former
  case, we have
  \begin{align*}
    G_n^B(e_B,b_2,\dots,b_n) &= H \circ
    m_2^A(e_A,G_{n-1}(b_2,\dots,b_n)) \\
    &= H \circ G_{n-1}(b_2,\dots,b_n) ,
  \end{align*}
  and in the latter case,
  \begin{align*}
    G_n^B(b_1,\dots,b_{n-1},e_B) &= H \circ
    m_2^A(G_{n-1}(b_1,\dots,b_{n-1}),e_A) \\
    &= (-1)^{|b_1|+\dots+|b_{n-1}|+1} H \circ
    G_{n-1}(b_1,\dots,b_{n-1}) .
  \end{align*}
  We see that $G_2(e_B,b)=(-1)^{|b|+1} G_2(b,e_B)=H\circ g(b)$,
  while if $n>2$, we see that $G_n(e_B,b_2,\dots,b_n)$ and
  $G_n(b_1,\dots,b_{n-1},e_B)$ lie in the image of $H\circ H$; in each
  case, these expressions vanish by the side conditions.
  
  A similar argument shows that $m_2(e_B,b)=(-1)^{|b|+1}
  m_2(b,e_B)=f\circ g(b)=b$, while if $n>2$, $m_n(e_B,b_2,\dots,b_n)$
  and $m_n(b_1,\dots,b_{n-1},e_B)$ lie in the image of $f\circ H$, and
  hence vanish by the side conditions.
\end{proof}

\subsection*{\Cinf-structures}
  
The analogues of commutative algebras in the setting of \Ainf-algebras
are \Cinf-algebras. These were introduced by Kadeishvili
\cite{Kadeishvili2}, who calls them commutative \Ainf-algebras; they
are called balanced \Ainf-algebras by Markl~\cite{Markl}. To define
\Cinf-structures, we recall the definition of the shuffle product on a
tensor algebra.

Let $A^\bull$ be a graded vector space, and let
\begin{equation*}
  T^nA = A^{\o n}
\end{equation*}
be the $n$th tensor power of $A^\bull$. Given elements
$(a_1,\dots,a_n)$ of $A$, we denote the element $a_1\o\dots\o a_n$ of
$T^nA$ by $a_1\dots a_n$, thinking of it as a word in the letters
$a_i$.  Let $TA$ be the (non-unital) tensor algebra
\begin{equation*}
  TA = \bigoplus_{n=1}^\infty T^nA .
\end{equation*}

The shuffle product is the graded commutative product on $TA$ defined
by the formula
\begin{equation*}
  a_1\dots a_p \shuffle a_{p+1}\dots a_{p+q} = \sum_{I\coprod J
    = \{1,\dots,n\} } (-1)^{\eps(I,J)} a(I,J) .
\end{equation*}
The word $a(I,J)$ equals $a_{\pi_1}\dots a_{\pi_{p+q}}$, where $\pi$
is the permutation defined in terms of $I = \{i_1<\dots <i_p\}$ and $J
= \{j_1<\dots <j_q\}$ by
\begin{equation*}
  \pi_k =
  \begin{cases}
    \ell , & i_\ell=k , \\
    p+\ell , & j_\ell=k ,
  \end{cases}
\end{equation*}
and the sign is given by
\begin{equation*}
  \eps(I,J)= \sum_{i_k>j_l} |a_{i_k}|\,|a_{p+j_l}| .
\end{equation*}
For example,
\begin{equation*}
  a_1a_2 \shuffle a_3 = a_1a_2a_3 +
  (-1)^{|a_2||a_3|} a_1a_3a_2 +
  (-1)^{(|a_1|+|a_2|)|a_3|} a_3a_1a_2 .
\end{equation*}

\begin{defn}
  A \Cinf-structure on a cochain complex $A^\bull$ is an
  \Ainf-structure such that for each $n>1$, $m_n^A$ vanishes on $TA
  \shuffle TA$.
\end{defn}

The basic example of a \Cinf-algebra is that of the \Ainf-algebra
associated to a differential graded commutative algebra, i.e.\ an
\Ainf-algebra $(A,m_{n})$ with $m_n=0$ for $n>2$, and such that $m_2$
satisfies
\begin{equation*}
  m_2(a,b) + (-1)^{|a|\,|b|} m_2(b,a) = 0 .
\end{equation*}
Inserting the definition \eqref{ab} of $ab$, we see that
\begin{equation*}
    (-1)^{|a|+1} ab + (-1)^{|a|\,|b|+|b|+1} ba = 0 ,
\end{equation*}
or equivalently,
\begin{equation*}
    ab = (-1)^{(|a|+1)(|b|+1)} ba .
\end{equation*}

\subsection*{The transfer of \Cinf-structures}

We now come to the main result of this article. When $A$ is a
differential graded commutative algebra and $B$ has vanishing
differential, this theorem was proved by Markl~\cite{Markl}.

\begin{thm}
  \label{main}
  If $A$ is a \Cinf-algebra, then the \Ainf-structure on $B$
  constructed by Kontsevich and Soibelman defines a \Cinf-structure.
\end{thm}

Let $\nabla_k:TA\to T^k(TA)$ be the morphism
\begin{equation*}
  a_1 \o \dots \o a_n \mapsto \sum_{\substack{n_1+\dots+n_k=n\\n_i>0}}
  (a_1\o\dots\o a_{n_1}) \o \dots \o ( a_{n-n_k+1}\o\dots\o a_n ) ,
\end{equation*}
and let $\mu:T(TA)\to TA$ be the morphism
\begin{equation*}
  (a_{11}\o\dots\o a_{1n_1})\o\dots\o(a_{k1}\o\dots\o a_{kn_k})
  \mapsto a_{11}\o\dots\o a_{kn_k} ,
\end{equation*}

We will use the following simple lemma in the proof of
Theorem~\ref{main}.
\begin{lem}
  \label{TT}
  We have
  \begin{multline*}
    \nabla_k(TA \shuffle TA ) \subset \\
    \mu ( T(TA) \shuffle T(TA) ) \oplus \bigoplus_{j=1}^k T^{j-1}(TA)
    \o (TA \shuffle TA) \o T^{k-j}(TA) .
  \end{multline*}
\end{lem}

\begin{proof}[Proof of Theorem \ref{main}]
  We will show, by induction on $n$, that $m_n^B(a\shuffle b)$ and
  $G_n(a\shuffle b)$ vanish on $T^nA \cap(TA \shuffle TA)$.
  
  By the induction hypothesis, the morphism $G_{n_1}\o\dots\o G_{n_k}$
  annihilates
  \begin{equation*}
    T^{j-1}(TA) \o (TA \shuffle TA) \o T^{n-j}TA .
  \end{equation*}
  Since it takes the subspace
  \begin{equation*}
    \mu ( T(TA) \shuffle T(TA) )
  \end{equation*}
  of $TA$ to $TA \shuffle TA$, it follows from Lemma \ref{TT} that
  \begin{equation*}
    (G_{n_1}\o\dots\o G_{n_k})(TA \shuffle TA ) \subset TA \shuffle TA
    .
  \end{equation*}
  By the induction hypothesis, this space is annihilated by $m_k^A$,
  and we conclude that $m_n^B(a\shuffle b)$ and $G_n(a\shuffle b)$
  vanish.
\end{proof}

\subsection*{Application to cosimplicial differential graded
  commutative algebras}

Let $\k[t_0,\dots,t_n,dt_0,\dots,dt_n]$ be the free graded commutative
algebra over $\k$ whose generators $t_i$ and $dt_i$ have degree zero
and $1$ respectively, and with differential $d$ defined on generators
by $d(t_i)=dt_i$ and $d(dt_i)=0$. Let $\Om_\bull$ be the simplicial
differential graded (dg) commutative algebra such that $\Om_n$ is the
dg algebra of differential forms with polynomial coefficient on the
$n$-simplex:
\begin{equation*}
  \Om_n = \frac{\k[t_0,\dots,t_n,dt_0,\dots,dt_n] } {
    ( t_0+\dots+t_n - 1 , dt_0 + \dots + dt_n ) } .
\end{equation*}

Let $A^\bull$ be a cosimplicial dg commutative algebra over $\k$.
\begin{exmp}
  If $X_\bull$ is a simplicial set, then $\k^{X_\bull}$ is the
  cosimplicial commutative algebra whose $n$-simplices are functions
  from $X_n$ to $\k$.
\end{exmp}

The Thom-Whitney normalization of $A^\bull$ (Navarro Aznar \cite{NA})
is the coend
\begin{equation*}
  \TotTW(A^\bull) = \tint^{n\in\Delta} \Om_n \o A^n . 
\end{equation*}
This cochain complex has a natural structure of a dg commutative
algebra. For example, $\TotTW(\k^{X_\bull})$ is the dg commutative
algebra $\Om(X_\bull)$ of polynomial coefficient differential forms on
$X_\bull$.

Let $\N_\bull$ be the simplicial cochain complex such that $\N_n$ is the
normalized simplicial cochain complex of the standard $n$-simplex
$\Delta^n$. The normalization of the cosimplical cochain complex
underlying the cosimplicial dg commutative algebra is the coend
\begin{equation*}
  \Tot(A^\bull) = \tint^{n\in\Delta} \N_n \o A^n . 
\end{equation*}
For example, $\Tot(\k^{X_\bull})$ is the complex $N(X_\bull)$ of
normalized simplicial cochains on $X_\bull$.

The integral $I_n:\Om_n^n\to\k$ over the $n$-simplex $\Delta^n$ is
given by the following explicit formula
\begin{equation*}
  I_n(t_1^{a_1}\dots t_n^{a_n}dt_1\dots dt_n)= \frac{a_1!\dots
    a_n!}{(a_1+\dots +a_n+n)!}.
\end{equation*}
Given a face $\{i_0<\dots<i_k\}$ of $\Delta^n$, let
\begin{equation*}
  I_{(i_{0}\dots i_{k})} : \Om_n^k \to \k
\end{equation*}
be the integral over the corresponding geometric $k$-simplex. Stokes's
theorem shows that
\begin{equation*}
  I_{(i_0\dots i_k)}(d\om) = \sum_{j=0}^k (-1)^j I_{(i_0\dots
    \widehat{\imath}_j\dots i_k)}(\om).
\end{equation*}
These maps give rise to a morphism of cochain complexes $f_n:\Om_n \to
\N_n$. Since these maps are compatible with the simplicial maps between
these complexes, we obtain a simplicial map $f_\bull:\Om_\bull \to
\N_\bull$.

The elementary differential forms (Whitney \cite{Whitney})
\begin{equation*}
  \om_{i_0\dots i_k} = k! \sum_{j=0}^k (-1)^j
  t_{i_j}dt_{i_0}\dots \widehat{dt}_{i_j}\dots dt_{i_k}
\end{equation*}
span a subcomplex of $\Om_n$ isomorphic to $\N_n$; this inclusion is
compatible with the simplicial maps between these complexes, and
defines a simplicial map $g_\bull:\N_\bull\hookrightarrow\Om_\bull$,
satisfying $f_\bull\circ g_\bull=1_{\N_\bull}$.

Given a cosimplicial dg commutative algebra $A^\bull$, these maps
induce quasi-isomorphisms
\begin{align*}
  f &: \TotTW(A^\bull) \to \Tot(A^\bull) , & g : \Tot(A^\bull) \to
  \TotTW(A^\bull) .
\end{align*}

Whitney shows (\cite{Whitney}, Chap.\ IV, Sect.\ 29) that the graded
commutative product $a\sqcup b=f(ga \wedge gb)$ on $\Tot(A^\bull)$
induces the cup product on cohomology. He does this by establishing
three conditions that characterize the cup product (Whitney,
\cite{Whitney:product}, Sect.\ 5):
\begin{enumerate}
\item the product $\alpha\sqcup\beta$ of two cochains is only nonzero
  on simplices which lie in the star of the supports of both $\alpha$
  and $\beta$;
\item if $\alpha$ is an $i$-cochain, then $\delta(\alpha\sqcup\beta) =
  \delta\alpha\sqcup\beta+(-1)^i\alpha\sqcup\delta\beta$;
\item the product $\sqcup$ has the $0$-cochain $1$ as its identity.
\end{enumerate}
(Strictly speaking, his proof is in the special case of cochains on a
simplicial complex, but it extends to the general case by naturality.)
However, while graded commutative, the product $\alpha\sqcup\beta$ is
not associative. We now show that it is part of a natural
\Cinf-structure on $\Tot(A^\bull)$, and in particular, is homotopy
associative.

Dupont constructed an explicit simplicial contraction of the dg
commutative algebra $\Om_\bull$ to $\N_\bull$, in this way giving a
rather concrete proof of de~Rham's theorem (Dupont
\cite{Dupont1,Dupont2}). We now recall his formula for this
contraction, although the explicit formula is not needed for the proof
of our main theorem.

For each vertex $e_i$ of the $n$-simplex $\Delta^n$, define the
dilation map
\begin{equation*}
  \phi_i : [0,1] \times \Delta^n \to \Delta^n
\end{equation*}
by the formula
\begin{equation*}
  \phi_i(u,t_0\dots t_n) = ((1-u)t_0,\dots,(1-u)t_i+u,\dots,(1-u)t_n) .
\end{equation*}
Let $(\pi_i){}_*:\Om^*([0,1]\times\Delta^n)\to\Om^{*-1}_n$ be
integration over the first factor. Define the operator
\begin{equation*}
  h_n^i : \Om_n^* \to \Om_n^{*-1}
\end{equation*}
by the formula
\begin{equation*}
  h_n^i\om = (\pi_i){}_\ast \phi_i^\ast \om .
\end{equation*}
Let $\eps_n^i:\Om_n\to\k$ be evaluation at the vertex $e_i$. The
Poincar\'e Lemma states that
\begin{equation*}
  \id_{\Om_n} - \eps_n^i =  dh_n^i + h_n^id .
\end{equation*}.

\begin{thm}[Dupont \cite{Dupont1,Dupont2}]
  The operator
  \begin{equation*}
    s_n = \sum_{k=0}^{n-1} \sum_{i_0<\dots <i_k} \om_{i_0\dots i_k}
    h^{i_k}\dots h^{i_0}
  \end{equation*}
  defines a simplicial contraction from $\Om_\bull$ to $\N_\bull$:
  \begin{equation*}
    \id_{\Om_n} - g_n\circ f_n = ds_n+s_nd .
  \end{equation*}
\end{thm}

\begin{thm}[Getzler \cite{Getzler}]
  \label{HH}
  The side conditions $f_n\circ s_n = 0$ and $s_n\circ s_n = 0$ hold.
\end{thm}

A cosimplicial dg commutative algebra $A^\bull$ is unital if each of
its component algebras $A^n$ is unital, and if the identities are
preserved by the cosimplicial maps.

\begin{thm}
  \label{cinf}
  The cochain complex $\Tot(A^\bull)$ associated to a cosimplicial
  commutative algebra $A^\bull$ (or rather, the shifted complex
  $\Tot(A^\bull)[1]$\,) has a canonical \Cinf-structure.

  If $A^\bull$ is unital, then $\Tot(A^\bull)[1]$ is unital.
\end{thm}
\begin{proof}
  The simplicial morphisms $f_\bull$, $g_\bull$ and $s_\bull$ induce a
  contraction $(f,g,H)$ from $\TotTW(A^\bull)$ to $\Tot(A^\bull)$ ---
  in particular, $H$ is the negative of $s$. The side conditions
  $f\circ H=0$ and $H\circ H=0$ hold by Theorem~\ref{HH}, and the side
  condition $H(1)=0$ holds, since $s_n(1)=0$ for all $n\ge0$.
  
  By Theorem \ref{main}, we may transfer the \Cinf-structure on
  $\TotTW(A^\bull)[1]$ induced by the dg commutative algebra structure
  on $\TotTW(A^\bull)$ to a canonical \Cinf-structure on
  $\Tot(A^\bull)$[1].
  
  By Theorem~\ref{unital}, the induced \Cinf-structure on
  $\Tot(A^\bull)[1]$ is unital if $A^\bull$ is unital.
\end{proof}

\begin{rem}
  Sullivan \cite{Sullivan} has given a different, less explicit,
  construction of this \Cinf-structure using obstruction theory.
  Wilson \cite{Wilson} has shown that in the limit where a simplicial
  complex $X_\bull$ is repeatedly subdivided, $m_2$ converges to the
  wedge product, and the higher homotopies converge to zero.
\end{rem}
  
We conclude by proving the conjecture of Remark A.3\,(3) of Sullivan
\cite{Sullivan}; the method of proof follows Fiorenza and Manetti
\cite{FM}. Mn\"ev~\cite{Mnev} has independently obtained this result
by the same method; he also studies the case of the circle. For an
entirely different approach, see Lawrence and Sullivan \cite{LS}.

By Theorem~\ref{cinf}, the complex of normalized simplicials cochains
on a simplicial set $X_\bull$ has a canonical unital \Cinf-structure.
In particular, taking $X=\Delta^1$, we obtain a \Cinf-structure on the
the complex $\N_1$ of normalized simplicial cochains on the interval.
  
Identify the cochain complex $\N_1$ with the subcomplex of $\Om_1$
spanned by the elementary differential forms $\{1,t,dt\}$. Since $1$
is an identity, it suffices to describe the products involving only
$t$ and $dt$.
\begin{prop}
  We have $m_2(t,t) = t$ and
  \begin{equation*}
    \textstyle
    m_{n+1}(dt^{\o i},t,dt^{\o n-i}) = (-1)^{n-i} \binom{n}{i} \,
    m_{n+1}(t,dt,\dots,dt) ,
  \end{equation*}
  where
  \begin{equation*}
    m_{n+1}(t,dt,\dots,dt) = \frac{B_n}{n!} \,dt .
  \end{equation*}
  All other products involving only $t$ and $dt$ vanish.
\end{prop}
\begin{proof}
  We have
  \begin{equation*}
    f(a)=a(0)\*1+(a(1)-a(0))\*t ,
  \end{equation*}
  hence $m_2(t,t)=f(t^2)=t$.      
    
  The product $m_{n+1}(t,dt,\dots,dt)$ may be calculated by induction
  on $n$, using the following explicit formula for Dupont's operator
  $s$:
  \begin{equation*}
    s ( t^k\,dt ) = \frac{t^{k+1} - t}{k+1} .
  \end{equation*}
  Let $p_n(t)$ be the sequence of polynomials defined by the
  recursion $p_1(t)=t$, and
  \begin{equation*}
    p_n(t) = s ( p_{n-1}(t) \, dt ) .
  \end{equation*}
  Then we have
  \begin{equation*}
    b_n = (-1)^{n-1} \int_0^1 p_n(t) \, dt .
  \end{equation*}
  But the recursion for $p_n(t)$ may be solved explicitly: we write
  \begin{equation*}
    P(z,t) = \sum_{n=0}^\infty z^n p_n(t) .
  \end{equation*}
  Then we have
  \begin{equation*}
    \frac{d}{dt} P(z,t) = z \Bigl( P(z,t) - \int_0^1 P(z,t) \, dt
    \Bigr) + z .
  \end{equation*}
  This equation is easily solved: we have
  \begin{equation*}
    P(z,t) = z \Bigl( \frac{e^{zt} - 1}{e^z - 1} \Bigr) ,
  \end{equation*}
  and hence $p_n(t) = ( B_n(t)-B_n )/n!$ and
  \begin{equation*}
    b_n = (-1)^n B_n/n! ,
  \end{equation*}
  where $B_n(t)$ is the $n$th Bernoulli polynomial, and $B_n=B_n(0)$
  is the $n$th Bernoulli number.
  
  The only planar trees contributing to the product
  \begin{equation*}
    m_{n+1}(dt^{\o i},t,dt^{\o n-i})
  \end{equation*}
  are binary planar trees with $n+1$ tails in which the path from the
  $(i+1)$st tail to the root passes through all $n$ vertices of the
  graph.  Such trees are enumerated by words of length $n$ in the
  letters $L$ and $R$, with $i$ instances of $R$, expressing whether
  at the $j$th vertex, the path from the $(i+1)$st tail to the root
  passes along the left ($L$) or right ($R$) branch. Each of these
  planar trees contributes the quantity
  $(-1)^i\,m_{n+1}(t,dt,\dots,dt)$ to the product.
\end{proof}

\end{document}